\documentclass[preprint,authoryear,natbib,fleqn]{article}
\usepackage{titlesec}

\usepackage{amsmath}
\usepackage{amssymb}
\usepackage{theorem}
\usepackage{xspace}
\usepackage{geometry}
\usepackage[utf8]{inputenc}    % pour les accents 
\def\pth#1{\left(#1\right)}

\def\cro#1{\left[#1\right]}

\DeclareMathAccent{\widehat}{\mathord}{largesymbols}{"62}

\def\ebo{\textrm{\mathversion{bold}$\mathbf{\beta^0}$\mathversion{normal}}}

\def\oo{\textrm{\mathversion{bold}$\mathbf{0}$\mathversion{normal}}}
\def\eb{\textrm{\mathversion{bold}$\mathbf{\beta}$\mathversion{normal}}}  
\def\ed{\textrm{\mathversion{bold}$\mathbf{\delta}$\mathversion{normal}}}

\def\eS{\textrm{\mathversion{bold}$\mathbf{\Sigma}$\mathversion{normal}}}

\def\eE{I\!\!E}

\def\e1{1\!\!1}

\def\eu{u'}

\def\XX{\textrm{\mathversion{bold}$\mathbf{X}$\mathversion{normal}}}

\def\xx{\textrm{\mathversion{bold}$\mathbf{x}$\mathversion{normal}}}
\def\hh{ \hspace*{0.5cm}}
\theoremstyle{plain}

\newtheorem{theorem}{Theorem}              %[section]
\newtheorem{lemma}{Lemma}               %[section]
\newtheorem{remark}{Remark}               %[section]
\newcommand{\beqn}{\begin{eqnarray*}}
\newcommand{\eeqn}{\end{eqnarray*}}
  
\def\ee1{\textrm{\mathversion{bold}$\mathbf{\varepsilon}$\mathversion{normal}}}  
\def\eth{\textrm{\mathversion{bold}$\mathbf{\theta}$\mathversion{normal}}}

\def\eu{\mathbf{{u}}}

\newcommand{\N}{\mathbb{N}}
\newcommand{\R}{\mathbb{R}}
\newcommand{\PP}{\mathbb{P}}

\def\argmin{\mathop{\mathrm{arg\,min}}} 
\date{}
\begin{document}
%--------------------------------------------------
%\begin{frontmatter}
\title {Model selection in high-dimensional quantile regression   with  seamless $L_0$ penalty}
\author{Gabriela CIUPERCA}
\maketitle

\begin{center}
{\it Universit\'e de Lyon, Universit\'e Lyon 1, 
CNRS, UMR 5208, Institut Camille Jordan, 
Bat.  Braconnier,\\
 43,  blvd du 11 novembre 1918, F - 69622 Villeurbanne Cedex, France,}
\end{center}
\begin{abstract}
In this paper we are interested in parameters estimation   of linear model when number of parameters increases with sample size. Without any assumption about moments of the model error, we propose and study the seamless $L_0$ quantile estimator. For this estimator we first give the convergence rate. Afterwards, we prove that it correctly distinguishes between zero and nonzero parameters and that the estimators of the nonzero parameters are asymptotically normal. A consistent BIC criterion to select the  tuning parameters is given. 
\end{abstract}
%%%%%%%%%%%%%%%
%%%%%%%%%%%%%%%
\textit{Keywords:}
High-dimension, quantile regression, seamless $L_0$ penalty, oracle properties, BIC criterion.  \\
\textit{MSC}: Primary 62F12; secondary  62J05.

%%%%%%%%%%%%%%%%%%%%%%
%%%%%%%%%%%%%%%%%%%%%%
%\maketitle
%\end{frontmatter}

\section{Introduction}
Consider a model where the number of regressors can increase with the sample size $n$:
\begin{equation}
\label{eq1}
Y_i=\XX_i^t \eb_n +\varepsilon_i, \qquad i=1, \cdots, n,
\end{equation}
where $\eb_n=(\beta_1, \cdots , \beta_{d_n})\in \R^{d_n}$ contains the regression parameters. The design vector $\XX_i$, for  observation $i$, is a deterministic vector of  dimension $d_n \times 1$. The random variable $\varepsilon_i$ is the model error. Denote by  $ \eb^0_n=(\beta^0_1, \cdots ,\beta ^0_{d_n})$ the true value, unknown,  of the parameter $\eb_n$. 
In order to automatically select the non-zero components of $\eb_n$ (therefore to select the significant variables), intuitively, the random optimization process  would penalize  with the "norm" $L_0$ (it is not a norm) defined by  $\|\eb_n\|_0=\sum^{d_n}_{j=1} \e1_{\beta_j \neq 0}$. This "norm"  has the disadvantage that it is not continuous in 0, then it  is computationally infeasible, since all possible models should be considered (all possible combinations of $\beta_j \neq 0$).
In this paper, we estimate the parameter  $\eb_n$ of (\ref{eq1}), penalizing the  quantile process with a   seamless  $L_0$ norm.   The difficulty in studing of  this type of estimation method is that  the  quantile process is  convex in $\eb_n$ and the seamless  $L_0$ penalty is concave.\\
\hh In literature on the high-dimension models, it was considered only the case of a quantile process penalized with a convex penalty  of type $L_1$.  Models with the number of variables exceeding the sample size ($d_n >n$) are studied by \cite{BelloniChernozhukov11}, \cite{FFB14}, \cite{ZGK13}. If $d_n <n$,  references \cite{WuLiu09}, \cite{ZouYuan08}  considered variable selection in a quantile model with convex penalties.\\
Penalized random process of type:
\begin{equation}
%\label{PP}
G_n(\eb_n)+{\cal P}en(\eb_n),
\end{equation}
with the process  $G_n(\eb_n)$ convex in $\eb_n$ and the penalty ${\cal P}en(\eb_n)$ nonconvex has been few studied.  
In \cite{FangPeng04}, $G_n(\eb)$ is $-\log$likelihood and the penalty is  nonconvex, with $d_n^5/n \rightarrow \infty$, as $n \rightarrow \infty$. 
For $d_n \gg n$, \cite{WLZ14} considered, for the particular case of  $Y|\XX=\xx$ sub-Gaussian,   $G_n(\eb_n)$ a loss function and ${\cal P}en(\eb_n)$ nonconvex loss penalty. For always $d_n \gg n$,  \cite{ZhangZhang12} considered $G_n(\eb_n)=(2n)^{-1} \sum^n_{i=1}(Y_i -\XX_i^t \eb_n)^2$ and  ${\cal P}en(\eb_n)$ concave. \cite{FXZ14} proposed an estimation method based on one-step local linear approximation, when the support set for $\eb^0_n$ is known.\\
\hh To overcome the disadvantage of the discontinuity in 0 of the norm $L_0$, \cite{DHL13} propose  a seamless  $L_0$ penalty:
\begin{equation}
\label{eq3}
{\cal P}en(\eb_n) \equiv \frac{\lambda_n }{ \log 2} \sum^{d_n}_{j=1} \log \big( \frac{ |\beta_j|}{|\beta_j| +\gamma_n} +1\big),
\end{equation}
with  $\lambda_n, \gamma_n >0$ two tuning parameters. If $\gamma_n \rightarrow 0$, the penalty   $L_0$ is obtained. Reference  \cite{DHL13} considers  $G_n(\eb_n)= n^{-1} \sum^n_{i=1}(Y_i -\XX_i^t \eb_n)^2$ , with  $(\varepsilon_i)$ i.i.d., $\eE[\varepsilon_i]=0$, $Var(\varepsilon_i)= \sigma^2$, suppositions  under  which the sparsity and the asymptotic normality of estimators are proved, if  $d_n/n \rightarrow 0$, for $n \rightarrow \infty$. If $Y$ belongs to the exponential family, \cite{LWL12} considers $G_n(\eb_n)= -  \log$likelihood$/n$, with  penalty  (\ref{eq3}), but with a stronger constraint on $d_n$: $d_n^5/n \rightarrow 0$ for $n \rightarrow \infty$. \\
\hh If the law of the error   $\varepsilon$ is unknown, or if the assumptions on the first two moments of the error are not satisfied, then the likelihood, least squares  methods with seamless $L_0$ penalty can not be used. This justifies the interest of the present paper, where  quantile process will be penalized with  seamless $L_0$ penalty (\ref{eq3}).\\
\hh We give some general notations. Throughout the paper, $C$ denotes a positives generic constant not dependent on $n$ which may take different values in different formula or even in different parts of the same formula.  All vectors and matrices are in bold and all vectors are column. For a vector $\mathbf{v}$, $\|\mathbf{v}\|_2$ is the Euclidean norm, $\mathbf{v}^t$ denotes the transposed of $\mathbf{v}$. For a matrix $\mathbf{M}$, $\|\mathbf{M}\|_2$ is the subordinate norm to the vector norm $\| .\|_2$, $\lambda_{\min}(\mathbf{M})$ and $\lambda_{\max}(\mathbf{M})$ are smallest and largest eigenvalues. We use also the notation $sgn(.)$ for the sign function and $tr(.)$ for the trace operator.\\
\hh The paper is organized as follows. In Section 2, we introduce and study the convergence rate, oracle properties of the seamless $L_0$ quantile estimator. In Section 3 we propose a consistent BIC criterion to select the tuning parameters. Finally, in Section 3, we present two lemmas useful to prove the main results.

\section{Seamless $L_0$ quantile estimator }
In this section we propose and study the seamless $L_0$ quantile estimator. 
For a fixed quantile index $\tau \in (0,1)$, the seamless $L_0$ quantile estimator is the parameter  which minimizes the process
\[
%\label{eq2}
Q_n (\eb_n) \equiv \frac{1}{2n} \sum^n_{i=1} \rho_\tau (Y_i -\XX_i^t \eb_n)+   \sum^{d_n}_{j=1} p_{SELO}(\beta_j),
\]
with the function $\rho_\tau(.): \R \rightarrow \R_{+}$ defined by  $\rho_\tau(u)= u(\tau -\e1_{u< 0})$ and for  $\beta \in \R$,
\[
p_{SELO}(\beta) \equiv \frac{\lambda_n}{\log 2} \log \big( \frac{|\beta|}{|\beta|+\gamma_n}+1\big).
\]
Then, the seamless $L_0$ quantile estimator is
\begin{equation}
\label{bb}
\widehat \eb_n \equiv \argmin_{\eb_n \in \R^{d_n}} Q_n(\eb_n).
\end{equation}
\begin{remark}
We emphasize that the results of \cite{FXZ14}, where a concave penalty is considered for quantile process, cannot be applied in the present paper, because our penalty cannot written as $\|\mathbf{c} \circ \eb\|_1$.
\end{remark}
For errors $(\varepsilon_i)$ of model (\ref{eq1}), we consider the following assumption:\\
\textit{(A1)}  $(\varepsilon_i)_{1 \leq i \leq n}$ are i.i.d., with  the distribution function  $F$ and density function $f$. The density function $f$ is continuously,  strictly positive in a neighborhood of zero and has a bounded first derivative in the neighborhood of 0. The $\tau$th quantile of $\varepsilon_i$ is zero: $\tau= F(0)$. \\
Let us denote $\alpha_n=(d_n/n)^{1/2}$. For the deterministic design $(\XX_i)_{1 \leq i \leq n}$ we suppose that:\\
\textit{(A2)} there exist constants $0 < r_0 \leq R_0 < \infty$ such that  $ r_0 \leq \lambda_{\min} ( {n}^{-1} \sum^n_{i=1} \XX_i \XX_i^t) \leq \lambda_{\max} ( {n}^{-1} \sum^n_{i=1} \XX_i \XX_i^t) \leq R_0 $. \\
\textit{(A3)}   $ \max_{1 \leq i \leq n} \| \XX_i\|_2 =o(\alpha_n^{-1})$.\\
On the  tuning parameters $\lambda_n$, $\gamma_n$ and on the  dimension $d_n$, we suppose:\\
\textit{(A4)} $d_n$ is such that ${\displaystyle  {d_n}/{n} \rightarrow 0}$, as $n \rightarrow \infty$.\\
\textit{(A5)} $\lambda_n =O(1)$, $\lambda_n \sqrt{{n}/{d_n}} \rightarrow \infty$ and $\gamma_n =O ({\displaystyle  d_n^{1/2}n^{-3/2}})$.\\

Assumptions (A1), (A2) are standard for linear model and (A3) is classic for an high-dimensional model. Assumptions (A4), (A5) are needed for statistical inference study of $\widehat{\eb}_n$ (see e.g. \cite{DHL13}, \cite{LHP14}).\\

For  $ \eb_n \in \R^{d_n}$, let be the difference between two  quantile processes:
\begin{equation}
\label{GnB}
G_n(  \eb_n) \equiv \sum^n_{i=1} [\rho_\tau(Y_i - \XX^t_i   \eb_n)-\rho_\tau(\varepsilon_i)].
\end{equation}
%Give also a relation true for all  $ a \in \R$, which will use in the result proofs:
%\[
%\label{eq4}
%\int^a_0 [\e1_{\varepsilon < t} - \e1_{\varepsilon < 0}] dt = (|a|-\varepsilon) \e1_{\min(0,a) \leq \varepsilon \leq \max(0,a)}.
%\]
Following theorem states that the estimators  $\widehat \eb_n$ has a convergence rate of order $\alpha_n$. If $d_n$ is bounded, we find the classic convergence rate $n^{-1/2}$ of quantile estimator for a finite-dimensional model (see \cite{Knight98}).

\begin{theorem}
\label{Lemma A.1}
Under assumptions (A1)-(A5), we have: $\| \widehat \eb_n -\eb^0_n \|_2 = O_{\PP}(\alpha_n)$.
\end{theorem}
\noindent {\bf Proof.} 
In order to prove the theorem, we show that for all $  \epsilon \in (0,1)$, there exists a constant large enough $ B>0$, such that we have, for any  $n$ large enough,
\begin{equation}
\label{t11}
\PP \cro{ \inf_{ \| \eu \|_2 =1} Q_n(\eb^0_n +B \alpha_n \eu) > Q_n(\eb^0_n)} \geq 1 - \epsilon .
\end{equation}
Fix $\epsilon \in (0,1)$ and consider some $ \eu  =(u_1, \cdots, u_{d_n})\in \R^{d_n}$,  with $\| \eu \|_2 =1$. Let be some constant $c>0$. Consider     
\begin{equation}
 Q_n(\eb^0_n +c \alpha_n \eu) -  Q_n(\eb^0_n) =  \frac{1}{2 n} G_n(c \alpha_n \eu) + \sum^{d_n}_{j=1} [p_{SELO}(\beta^0_j +c \alpha_n u_j)-p_{SELO}(\beta^0_j)].
\label{eq5}
\end{equation}
For the   penalty, we have the following inequality:
\[
\sum^{d_n}_{j=1} [p_{SELO}(\beta^0_j+c \alpha_n u_j)- p_{SELO}(\beta^0_j)] \geq \sum_{j \in J(\eu)} [p_{SELO}(\beta^0_j+c \alpha_n u_j)- p_{SELO}(\beta^0_j)] ,
\]
where $J(\eu) \equiv \{l \in \{ 1, \cdots, d_n \} ; \; p_{SELO}(\beta^0_l+c \alpha_n u_l)- p_{SELO}(\beta^0_l) <0 \}$. 
Because  $c$, $\eu$ are fixed and $\alpha_n \rightarrow 0$ , then by  Lemma \ref{Lemma gh}, for all  $j \in J(\eu)$, and for large enough $n$,  there exists $\tilde C_j >0$ such that 
\[
p_{SELO}(\beta^0_j+c \alpha_n u_j)- p_{SELO}(\beta^0_j) = \frac{\lambda_n}{\log 2} \cro{ g(\beta^0_j+c \alpha_n u_j)-g(\beta^0_j)}=\frac{\lambda_n}{\log 2}  \tilde C_j \alpha_n |u_j| \gamma_n (-1)^{sgn(\beta^0_j(\beta^0_j+ c  \alpha_n u_j) )}.
\]
Thus, by assumptions (A4) and (A5), we have:
\begin{equation}
\label{pl}
%\begin{array}{lll}
  \sum_{j \in J(\eu)} [p_{SELO}(\beta^0_j+c \alpha_n u_j)- p_{SELO}(\beta^0_j) ] >   - \frac{\lambda_n \alpha_n \gamma_n }{\log 2}\sum_{j \in J(\eu)} \tilde{C_j} |u_j|  
 = -O(\lambda_n \alpha_n \gamma_n d_n) = - O(\alpha_n \alpha_n^{3/2}) = -o(\alpha_n^2).
%\end{array}
\end{equation}
We now study the expectation of  $G_n(c \alpha_n \eu)$:
% \begin{eqnarray}
% \displaystyle {\frac{1}{n} \eE[G_n(C\alpha_n \eu)] } &=& \displaystyle {\frac{1}{n} \sum^n_{i=1}  \eE[ \rho_\tau(\varepsilon_i-C\alpha_n  \XX_i^t\eu)-\rho_\tau(\varepsilon_i)]}      \nonumber \\
%   &=& \displaystyle { \frac{1}{n} \sum^n_{i=1} \eE\cro{ \int^{C\alpha_n \XX_i^t \eu}_0 \e1_{0 < \varepsilon_i <t} dt} }  \nonumber \\
 %  & = & \displaystyle { \frac{1}{n} \sum^n_{i=1} \int^{C\alpha_n \XX_i^t \eu}_0  [F(t)-F(0)] dt} .  \nonumber
%\end{eqnarray}
\[
 \eE[G_n(c\alpha_n \eu)] = \sum^n_{i=1}  \eE [ \rho_\tau(\varepsilon_i-c\alpha_n  \XX_i^t\eu)-\rho_\tau(\varepsilon_i)] =  \sum^n_{i=1} \eE \big[ \int^{c\alpha_n \XX_i^t \eu}_0 \e1_{0 < \varepsilon_i <t} dt \big] = \sum^n_{i=1} \int^{c\alpha_n \XX_i^t \eu}_0  [F(t)-F(0)] dt.
\]
On the other hand, by (A1),   for $v \rightarrow 0$, we have
$
\int^v_0[F(t)-F(0)]dt=\frac{f(0)}{2} v^2+o(v^2)$. 
Using  (A3), we have:
\[
 \frac{1}{n} \sum^n_{i=1} \int^{c \alpha_n \XX_i^t \eu}_0[F(t)-F(0)]dt= \frac{f(0)}{2} c  \alpha^2_n  \frac{1}{n} \sum^n_{i=1} (\XX_i^t \eu)^2 +o(\alpha^2_n  \frac{1}{n} \sum^n_{i=1} \eu^t (\XX_i \XX_i^t)\eu).
\]
Then
\begin{equation}
\label{ec}
\frac{1}{n} \eE[G_n(c \alpha_n \eu)] = c \frac{f(0)}{2} \alpha^2_n \frac{1}{n}  \eu^t (\sum^n_{i=1} \XX_i \XX_i^t) \eu (1+o(1)).
\end{equation}
Consider now the random  variables  ${\cal D}_i \equiv (1-\tau) \e1_{\varepsilon_i >0} - \tau \e1_{\varepsilon_i >0}$, 
$
R_i \equiv \rho_\tau(\varepsilon_i-c\alpha_n \XX_i^t \eu)-\rho_\tau(\varepsilon_i) - c \alpha_n  D_i \XX_i^t \eu
$
and the random vector
$
\textbf{W}_n \equiv \sum^n_{i=1} c \alpha_n D_i \XX_i^t$.  
Thus, the process $G_n$ can be written:
\[
G_n(c\alpha_n \eu) =\eE[G_n(c\alpha_n \eu)] +\textbf{W}_n \eu +\sum^n_{i=1}  \big[R_i-\eE[R_i] \big].
\]
But, since, by (A1), the errors   $(\varepsilon_i)$ are independent, using also $|R_i| \leq | c \alpha_n \XX_i^t \eu| \e1_{|\varepsilon_i| \leq | c \alpha_n \XX_i^t \eu|}$, we have:
 %\begin{eqnarray}
 %\displaystyle {\eE[\sum^n_{i=1}  [R_i-\eE[R_i]]^2 }&= & \displaystyle { \sum^n_{i=1}  \eE[R_i-\eE[R_i] ]^2}  \nonumber \\
 %& \leq & \displaystyle { \sum^n_{i=1} \eE[R_i]^2  }  \nonumber \\
 %& \leq & \displaystyle { \sum^n_{i=1} | C \alpha_n \XX_i^t \eu|^2 \eE[\e1_{|\varepsilon_i| \leq | C \alpha_n \XX_i^t \eu|}].} \nonumber
 %\end{eqnarray}
 \[
 \eE[\sum^n_{i=1}  [R_i-\eE[R_i]]^2 =  \sum^n_{i=1}  \eE[R_i-\eE[R_i] ]^2 \leq  \sum^n_{i=1} \eE[R_i]^2  \leq \sum^n_{i=1} | c \alpha_n \XX_i^t \eu|^2 \eE[\e1_{|\varepsilon_i| \leq | c \alpha_n \XX_i^t \eu|}].
 \]
Taking into account assumptions (A1) and  (A3), we have
 $
 \eE[\e1_{|\varepsilon_i| \leq | c \alpha_n \XX_i^t \eu|}] =F(| c \alpha_n \XX_i^t \eu|)- F(-| c \alpha_n \XX_i^t \eu|)= C \alpha_n |\XX_i^t \eu| \leq   C \alpha_n \max_{1 \leq i \leq n } \| \XX_i\|_2 =o(1)$, with $C>0$.
Then, using assumption (A2), we obtain: 
 \begin{equation}
 \label{er}
 \eE\big[\sum^n_{i=1}  [R_i-\eE[R_i]]\big]^2 =o \big( \alpha^2_n \eu^t \sum^n_{i=1} \XX_i \XX_i^t \eu \big) =o(d_n).
 \end{equation}
Consider now the random  variable  $U_n \equiv d_n^{-1/2}\sum^n_{i=1}  [R_i-\eE[R_i] ]$. Taking into account (\ref{er}), we have 
 $
 \eE[U^2_n]=o(1)
 $. Since $E[U_n]=0$, by  Bienaymé-Tchebychev inequality, we have $U_n \overset{\PP} {\underset{n \rightarrow \infty}{\longrightarrow}}  0$. Thus $\sum^n_{i=1}  \big[R_i-\eE[R_i]\big]=o_{\PP}(d_n^{1/2}) $. Returning to $G_n$,  we have, taking into account (\ref{ec}):
\[
G_n(c \alpha_n \eu)=  \eE[G_n(c \alpha_n \eu)]+\textbf{W}_n \eu+o_{\PP}(d_n^{1/2})
\]
or again
\begin{equation}
\label{Gn}
G_n(c \alpha_n \eu)= \bigg(\frac{f(0)}{2} c^2 d_n \eu^t (\frac{1}{n}\sum^n_{i=1} \XX_i \XX_i^t) \eu+d_n^{1/2} c \big( \sum^n_{i=1} \frac{{\cal D}_i \XX_i^t}{\sqrt{n}}\big) \eu \bigg)(1+ o_{\PP}(1))  +o_{\PP}(d_n^{1/2}).
\end{equation}
Since $n^{-1/2}\sum^n_{i=1}  {\cal D}_i \XX_i^t \eu$ converges in distribution to a centered normal distribution, by assumptions (A4) and (A2), for a large enough constant $B$, we have that the first term of the right side that will dominate in  (\ref{Gn}). Then,  
\begin{equation}
\label{ggg}
\frac{1}{2 n} G_n(B \alpha_n \eu)=f(0) B^2 \alpha^2_n \eu^t (\frac{1}{n}\sum^n_{i=1} \XX_i \XX_i^t) \eu (1+ o_{\PP}(1)) .
\end{equation}
Thus, for $n$ and  $B$ large enough, we have  $(2 n)^{-1} G_n(B \alpha_n \eu) >0$. 
On the other hand, by relations (\ref{eq5}) and  (\ref{pl}), 
$
 Q_n(\eb^0_n +B \alpha_n \eu) -  Q_n(\eb^0_n) >  (2 n)^{-1} G_n(B \alpha_n \eu) -o(\alpha_n^2). 
$
Taking also into account relation (\ref{ggg}) and assumption (A2), we obtain (\ref{t11}).
\hspace*{\fill}$\blacksquare$ \\

Let us consider the parameter set, with the  constant  $B>0$ of relation (\ref{t11}):
\[
{\cal V}_{\alpha_n}(\eb^0_n) \equiv \{ \eb_n \in  \R^{d_n} ; \quad \| \eb_n - \eb^0_n \| \leq B \alpha_n \}.
\]
According to Theorem \ref{Lemma A.1}, the seamless $L_0$ quantile estimators belong to ${\cal V}_{\alpha_n}(\eb^0_n)$, with a probability converging to 1. \\
For the index set $\mathcal{A}$, with $\mathcal{A} \subseteq \{1, \cdots , d_n \}$, we will denote by $|\mathcal{A}|$ its cardinal. Throughout the paper, we denote by $\eb_{\mathcal{A}}$ the sub-vector of $\eb_n$ containing the corresponding components of $\mathcal{A}$. Similarly for $\XX_{i , {{\cal A}}}$. 
Consider also the following index set:
\begin{equation}
\label{A0}
{\cal A}^0 \equiv \{ j \in \{ 1, \cdots, d_n\} ; \quad  \beta^0_j \neq 0\} .
\end{equation}
The following theorem gives the  oracle properties for the estimators  $\widehat \eb_n=(\widehat \beta_1, \cdots, \widehat \beta_{d_n})$, defined by (\ref{bb}). Note that, with respect to the  paper of \cite{DHL13}, for showing the normality of the nonzero estimators,  the condition $\eE[|\varepsilon_i|^{2+\delta}] <M$ is not needed, for some $\delta >0$ and $M < \infty$.

\begin{theorem}
\label{theorem 1}
Under assumptions (A1)-(A5), we have \\
(i) $\displaystyle{ \lim_{n \rightarrow \infty} \PP[\{j \in \{ 1, \cdots, d_n\}; \; \widehat \beta_j \neq 0\} ={\cal A}^0]=1}$.\\
(ii) For any vector $\eu$ of dimension $|{\cal A}^0|$ such that $\|\eu\|_2=1$, if we denote $ \eS_{{\cal A}^0} \equiv n^{-1} \sum^n_{i=1} \XX_{i , {{\cal A}^0}} \XX^t_{i, {{\cal A}^0}}$, then  
 $$\sqrt{n}(\eu^t \Sigma^{-1}_{{\cal A}^0})^{-1} \eu)^{-1/2} \eu^t (\widehat \eb_{{\cal A}^0}- \eb^0_{{\cal A}^0}) \overset{\cal L} {\underset{n \rightarrow \infty}{\longrightarrow}} {\cal N}\big(0, \frac{\tau(1-\tau)}{f^2(0)}\big). $$
\end{theorem}
\noindent {\bf Proof}. 
\textit{(i)} If we denote by  ${{\cal A}^0}^c$ the complementary set of $\mathcal{A}^0$ in $\{ 1, \cdots , d_n \}$, we will prove that for any $\eb_n=(\eb_{{\cal A}^0},\eb_{{{\cal A}^0}^c}) \in {\cal V}_{\alpha_n}(\eb^0_n)$ such that $\|\eb_{{\cal A}^0} - \eb^0_{{\cal A}^0} \|_2=O_{\PP}(\alpha_n)$ and any constant $C>0$, we have
\begin{equation}
\label{Lemma A.2}
Q_n((\eb_{{\cal A}^0},\textbf{0}))=\min_{\|  \eb_{{{\cal A}^0}^c}\| \leq C \alpha_n} Q_n((\eb_{{\cal A}^0}, \eb_{{{\cal A}^0}^c})).
\end{equation}
Consider the following parameter set 
$
{\cal W}_n \equiv \{ \eb_n \in {\cal V}_{\alpha_n}(\eb^0_n);  \| \eb_{{{\cal A}^0}^c}\|_2>0 \}$.
We show that $\PP[\widehat \eb_n \in {\cal W}_n] \rightarrow 0$, as $n \rightarrow \infty$. \\
Let $\eb_n =(\eb_{{\cal A}^0}, \eb_{{{\cal A}^0}^c}) \in {\cal W}_n$  and an another parameter $\widetilde \eb_n =(\widetilde \eb_{{\cal A}^0}, \widetilde  \eb_{{{\cal A}^0}^c}) \in {\cal V}_{\alpha_n}(\eb^0_n)$, such that $\widetilde \eb_{{\cal A}^0}= \eb_{{\cal A}^0}$  and $\widetilde  \eb_{{{\cal A}^0}^c} = \textbf{0}$. 
Define
\begin{equation}
\label{DDn}
D_n(\eb_n, \widetilde \eb_n) \equiv  Q_n(\eb_n) - Q_n(\widetilde{\eb}_n) =\frac{1}{2 n} \sum^n_{i=1} \cro{\rho_\tau(Y_i- \XX^t_i \eb_n)-\rho_\tau(Y_i - \XX^t_i \widetilde \eb_n)} +\sum_{j \in {{\cal A}^0}^c} p_{SELO}(\beta_j).
\end{equation}
 Concerning the penalty of relation (\ref{DDn}), as in the proof of Lemma A.2 of \cite{DHL13}, relation (A.7), we have that there  exists $C >0$ such that
\[
\sum_{j  \in {{\cal A}^0}^c}p_{SELO}(\beta_j) \geq \frac{\lambda_n}{\log 2}  \log \big(\frac{C}{C+\gamma_n \alpha_n} +1 \big)  \| \eb_n- \widetilde \eb_n\|_2 .
\]
On the other hand, by assumption (A3), we have that there exists  $ C_1 >0$ such that
$
\liminf_{n \rightarrow \infty} \big( \log \big({C}/{(C+\gamma_n \alpha_n)} +1 \big)\big)> C_1 >0$.
Then,  for $n$ large enough, there exists $ \tilde C >0$ such that
\begin{equation}
\label{eq7}
\sum_{j  \in {{\cal A}^0}^c} \frac{p_{SELO}(\beta_j)}{\| \eb_n- \widetilde \eb_n\|_2} \geq \tilde{C} \lambda_n .
\end{equation}
Let be the identity that follows from \cite{Knight98}, for any $x, y \in \R$,  
\[
% \label{Knight}
\rho_\tau(x-y)- \rho_\tau(x)=y(\e1_{x \leq 0} - \tau)+\int^y_0 (\e1_{x \leq t} -\e1_{x \leq 0})dt.
\]
Using this  relation for the first sum of (\ref{DDn}), we obtain:
\begin{equation}
\begin{array}{c}
\displaystyle{ \frac{1}{2 n} \sum^n_{i=1} \big[\rho_\tau(Y_i- \XX^t_i \eb_n)-\rho_\tau(Y_i - \XX^t_i \widetilde \eb_n)\big]= \frac{1}{2 n} (\eb_n -\widetilde \eb_n)^t \sum^n_{i=1} \XX_i  [\e1_{Y_i - \XX^t_i \widetilde \eb_n \leq 0} - \tau]}  \\
\displaystyle  {+\frac{1}{n} \sum^n_{i=1} \int^{\XX_i^t (\eb_n - \widetilde \eb_n)}_0 [\e1_{Y_i - \XX^t_i \widetilde \eb_n \leq t} -\e1_{Y_i - \XX^t_i \widetilde \eb_n \leq 0}] dt} 
\equiv T_{1n}+T_{2n} 
\end{array} 
\label{eqq0}
\end{equation}
For $T_{1n}$ we have, by assumption (A3) and since the density $f$ is bounded in a neighborhood of  0:
\[
\eE[T_{1n}]= (\eb_n- \widetilde \eb_n)^t \frac{1}{2 n} \sum^n_{i=1} \XX_i [F(\XX^t_i (\widetilde \eb_n -\eb^0_n))-F(0)] = (\eb_n- \widetilde \eb_n)^t \frac{1}{2 n} \big( \sum^n_{i=1}\XX_i  \XX^t_i\big) (\eb^0_n - \widetilde \eb_n) f(0) (1+o(1)).   
\]
Then
$
 |\eE[T_{1n}]| \leq \|\eb_n- \widetilde \eb_n \|_2  \big\| (2n)^{-1}\sum^n_{i=1}\XX_i  \XX^t_i  \big\|_2 \| \eb^0_n - \widetilde \eb_n\|_2 f(0) (1+o(1)) $.
 Since the matrix  $ n^{-1}\sum^n_{i=1}\XX_i  \XX^t_i $ is Hermitian, we have that $\left\| n^{-1}\sum^n_{i=1}\XX_i  \XX^t_i  \right\|_2 = \lambda_{\max}(n^{-1}\sum^n_{i=1}\XX_i  \XX^t_i ) \leq R_0$. Hence, by (A2), we have 
$
|\eE[T_{1n}]| \leq \|\eb_n- \widetilde \eb_n \|_2  \| \eb^0_n - \widetilde \eb_n\|_2 R_0 f(0)$.
Therefore
$
\eE[T_{1n}]= O( \|\eb_n- \widetilde \eb_n \|_2 \| \eb^0_n - \widetilde \eb_n\|_2 )=O( \|\eb_n- \widetilde \eb_n \|_2^2 )$. 
By calculations analogous to $\eE[T_{1n}]$, using independence of $\varepsilon_i$, we have that 
$\eE[T^2_{1n}] = C n^{-1} {\|\eb_n- \widetilde \eb_n\|^3} \rightarrow 0$,  for $n \rightarrow \infty$. Since $Var[T_{1n}] \leq \eE[T^2_{1n}]$, using Bienaymé-Tchebychev inequality, we obtain
\begin{equation}
\label{eqq1}
T_{1n}= C \|\eb_n- \widetilde \eb_n\|^2_2 (1+o_{\PP}(1)).
\end{equation}
Study now   $T_{2n}$ of (\ref{eqq0}), which can be written as:
$
T_{2n}= n^{-1} \sum^n_{i=1} \int^{\XX^t_i (\eb_n - \widetilde \eb_n)}_0 [\e1_{\varepsilon_i \leq t- \XX_i^t (\eb^0_n - \widetilde \eb_n)} - \e1_{\varepsilon_i \leq - \XX_i^t (\eb^0_n - \widetilde \eb_n)} ] dt $.
Then, taking into account that $\eb_n \in {\cal V}_{\alpha_n}(\eb^0_n) $, together with assumptions (A1),  (A3), we have 
\[
\eE[T_{2n}]  =  \frac{1}{n} \sum^n_{i=1} \int^{\XX^t_i (\eb_n - \widetilde \eb_n)}_0[F(t- \XX_i^t (\eb^0_n - \widetilde \eb_n)) - F(- \XX_i^t (\eb^0_n - \widetilde \eb_n))] dt  
  =  \frac{1}{n} \sum^n_{i=1} \int^{\XX^t_i (\eb_n - \widetilde \eb_n)}_0 [t \cdot f(\XX_i^t (\widetilde \eb_n - \eb^0_n)) + o(t) ]dt .   
 \]
 By Theorem \ref{Lemma A.1}, together with assumptions (A1),   (A3), we have that $f(\XX_i^t (\widetilde \eb_n - \eb^0_n))$ is bounded by a constant $C \in (0,\infty)$. 
 Thus, as for $T_{1n}$, using assumption (A2) and the fact that 
 $
n^{-1} \sum^n_{i=1} \| \XX_i\|^2_2 -tr\big( n^{-1}  \sum^n_{i=1} \XX_i \XX_i^t\big) \rightarrow 0, 
 $
 we have
 $
 |\eE[T_{2n}]| \leq 
  {C}n^{-1} \sum^n_{i=1} \| \XX_i \|^2_2   \| \eb_n -\widetilde \eb_n \|_2  \| \widetilde \eb_n - \eb^0_n\|_2   +o\big( n^{-1} \sum^n_{i=1}\XX_i^t (\eb_n-\widetilde \eb_n)\big)  
 = C \| \eb_n -\widetilde \eb_n \|^2_2  $.
We show similarly that $\displaystyle{\eE[T_{2n}^2] = C n^{-1} \| \eb_n - \widetilde \eb_n\|^3}$. 
Then, by Bienaymé-Tchebychev inequality, we have:
\begin{equation}
\label{eqq2}
T_{2n}= C \| \eb_n - \widetilde \eb_n\|^2_2 (1+o_{\PP}(1)).
\end{equation}
Hence, by relations  (\ref{eqq1}), (\ref{eqq2}), we obtain
\[
T_{1n}+T_{2n} = C \| \eb_n - \widetilde \eb_n\|^2_2 (1+o_{\PP}(1)).
\]
Thus, taking into account this last relation together with relations (\ref{DDn}), (\ref{eq7}), (\ref{eqq0}), and since   $\eb_n \in {\cal W}_n$, we have: 
$
 {D_n(\eb_n, \widetilde \eb_n)}{\| \eb_n - \widetilde \eb_n \|_2^{-1}} \geq C \| \eb_n - \widetilde \eb_n\|_2+\widetilde{C} \lambda_n $.
Since $\| \eb_n - \widetilde \eb_n \| =O( \alpha_n)$ and $\lambda_n/\alpha_n \rightarrow \infty$ by (A5),   we have that there exists $C_+ >0$ such that
$
 {D_n(\eb_n, \widetilde \eb_n)}{\| \eb_n - \widetilde \eb_n \|^{-1}_2} > C_+ \lambda_n >0 $.
But for $\eb^0_n$, taking into account the definition of  $\widetilde \eb_n$ and that of $D_n(\eb^0_n, \widetilde \eb_n) \equiv Q_n(\eb^0_n) - Q_n(\widetilde \eb_n)$, we have that $D_n(\eb_n^0, \widetilde \eb_n)=C \| \eb_n^0 - \widetilde \eb_n\|^2_2(1+o_{\PP}(1))$. Then, by (A3), we have  $\PP[\eb_n \in {\cal W}_n] \rightarrow 0$ and  relation (\ref{Lemma A.2}) follows.\\

\textit{(ii)} Taking into account  the estimator  convergence rate  obtained by Theorem  \ref{Lemma A.1} and claim \textit{(i)}, the estimator $\widehat \eb_n$ can be written $\widehat \eb_n=\eb^0_n +\alpha_n \ed$, with,   $\ed \equiv(\delta_1, \cdots, \delta_{d_n}) \in \R^{d_n}$,  $\ed_{{{\cal A}^0}^c}=\oo$ and $\| \ed_{{\cal A}^0} \|^2_2 \leq C  |{{\cal A}^0}|$. Consider then
\begin{equation}
\label{eqq4}
Q_n(\eb^0_n +\alpha_n \ed)-Q_n(\eb^0_n) = \frac{1}{2 n} \sum^n_{i=1} [\rho_\tau(Y_i -\XX_i^t( \eb^0_n +\alpha_n \ed)) - \rho_\tau(\varepsilon_i)] +{\cal P}, 
\end{equation}
with ${\cal P} \equiv \sum_{j \in {\cal A}^0} p_{SELO}(\beta^0_j +\alpha_n \delta_j)- p_{SELO}(\beta^0_j)$. 
Let us first study ${\cal P}$.
For any  $  j \in {{\cal A}^0}$, by Lemma \ref{Lemma gh}, we have that there exists a constant   $\widetilde C_j$ such that  
\[ 
p_{SELO}(\beta^0_j +\alpha_n \delta_j)- p_{SELO}(\beta^0_j)=  \frac{\lambda_n} {\log 2} [g(\beta^0_j +\alpha_n \delta_j)-g(\beta^0_j)] = \frac{\lambda_n} {\log 2} \gamma_n  \pth{|\beta^0_j +\alpha_n \delta_j| - |\beta^0_j| }  \widetilde C_j, 
\]
with $|\widetilde C_j| < \infty$, for any $ j \in {{\cal A}^0}$. Since $\alpha_n  \rightarrow 0$, $|\beta^0_j| >C>0$, $\forall j \in {{\cal A}^0}$ and $\delta_j  $   bounded, we have that for $n$ large enough, the parameters  $ \beta^0_j +\alpha_n \delta_j$ and $\beta^0_j$ have the same sign. 
Then
\begin{equation}
\label{eqq3}
{\cal P} = C    \frac{\lambda_n} {\log 2} \gamma_n \alpha_n \sum_{j \in {{\cal A}^0}} (\pm \delta_j) = C \lambda_n \alpha_n \gamma_n |{{\cal A}^0}|.
\end{equation}
For the first term of the right-hand side of  (\ref{eqq4}) we have:
\[
\frac{1}{2n} \sum^n_{i =1} [\rho_\tau(\varepsilon_i - \alpha_n \XX^t_i \ed)-\rho_\tau(\varepsilon_i)]= \frac{\alpha_n}{2n} \sum^n_{i=1}  \XX^t_i \ed [ \e1_{\varepsilon_i \leq  0}-\tau] + \frac{1}{2n} \sum^n_{i=1} \int^{\alpha_n \XX^t_i \ed} [\e1_{\varepsilon_i \leq t}- \e1_{\varepsilon_i \leq 0}] dt \equiv J_1 +J_2.
\]
Since $\eE[J_1]=0$, using independence of $(\varepsilon_i)$, assumption (A4)  and the Cauchy-Schwarz inequality, we get that
\[
Var(J_1)  \leq  \eE[J^2_1]  =   \frac{\alpha^2_n}{4 n^2} \tau(1-\tau) \sum^n_{i=1} ( \XX^t_i \ed )^2   \leq   \frac{\alpha^2_n}{4 n^2} \tau(1-\tau) \sum^n_{i=1} \| \XX^t_{i {{\cal A}^0}}\|^2_2   \| \ed_{{\cal A}^0} \|_2^2 = \frac{C}{n}  \alpha^2_n |{{\cal A}^0}| \leq \alpha^2_n \frac{d_n}{n} \rightarrow 0. \]
 For $J_2$ we have:
\begin{equation}
\eE[J_2] = \frac{1}{2n} \sum^n_{i=1} \int^{\alpha_n \XX^t_i \ed}_0 \pth{t f(0) + o(t^2)  } dt  =\frac{1}{4} f(0)\alpha_n^2  \ed^t \big(\frac{1}{n} \sum^n_{i=1}\XX_i \XX^t_i  \big) \ed  (1+o(1)) .
 \label{eq16}
\end{equation}
Using assumption (A2), we have that
$
f(0)\alpha_n^2  \|\ed\|^2_2 \cdot \lambda_{\min}\big( n^{-1} \sum^n_{i=1}\XX_{i  } \XX^t_{i  }   \big) \leq \eE[J_2] \leq f(0)\alpha_n^2  \|\ed\|^2_2 \cdot \lambda_{\max}\big(n^{-1} \sum^n_{i=1}\XX_{i } \XX^t_{i  }   \big)$. 
Taking into account the fact that $\|\ed\|^2_2=\|\ed_{{\cal A}^0}\|^2_2 \leq C  |{{\cal A}^0}|$, we have $\eE[J_2]=C f(0)\alpha_n^2 |{{\cal A}^0}| $.  We prove similarly $Var(J_2)=O(n^{-1}\alpha^3_n | \mathcal{A} ^0| )$. 
We compare $\alpha^2_n  |{{\cal A}^0}| $ with $\lambda_n \alpha_n \gamma_n |{{\cal A}^0}|$ obtained by (\ref{eqq3}) for the penalty,
$
\frac{\alpha^2_n |{{\cal A}^0}|}{\lambda_n \alpha_n \gamma_n |{{\cal A}^0}|} =\frac{\alpha_n}{\lambda_n \gamma_n }$. 
By (A5), 
$ { \gamma_n = O\pth{ \frac{\alpha_n^3}{d_n}}}$, thus $ {\frac{\gamma_n}{\alpha_n }=\frac{\alpha^2_n}{d_n}} $. Then $ { \frac{\alpha_n}{\lambda_n \gamma_n } =\frac{d_n}{\alpha^2_n} \frac{1}{\lambda_n}= \frac{n}{\lambda_n} \rightarrow \infty}$, as $n \rightarrow \infty$.  Thus, minimizing (\ref{eqq4}) amounts to  minimizing $J_1+J_2$, with respect to $\alpha_n \ed$. Using (\ref{eq16}), we obtain:
\begin{equation}
\label{MM}
\frac{1}{2 n} \sum^n_{i=1} [\rho_\tau(Y_i -\XX_i^t( \eb^0_n +\alpha_n \ed)) - \rho_\tau(\varepsilon_i)] = \frac{\alpha_n}{2 n}  \sum^n_{i=1}\XX^t_{i,{{\cal A}^0}} \ed_{{\cal A}^0} [\e1_{\varepsilon_i <0} - \tau]+ \frac{1}{4} f(0) \alpha^2_n \ed^t_{{\cal A}^0} \eS_{{\cal A}^0} \ed_{{\cal A}^0}(1+o_{\PP}(1)).
\end{equation}
The minimizer of  (\ref{MM}) is:
\begin{equation}
\label{ad}
\alpha_n \ed_{{\cal A}^0}=-\frac{1}{n} \frac{1}{f(0)} \eS^{-1}_{{\cal A}^0} \big( \sum^n_{i=1} \XX_{i,{ {\cal A}^0}} (\e1_{\varepsilon_i \leq 0}-\tau)\big).
\end{equation}
For studying (\ref{ad}), let us consider the following   independent random variable sequence
$
W_i \equiv (f(0))^{-1} \eu^t \eS^{-1}_{\cal A}\XX_{i {\cal A}} (\e1_{\varepsilon_i \leq 0}-\tau)$, 
with $\eu $ a  vector of  dimension $|{{\cal A}^0}|$ and such that  $\|\eu\|_2=1$. We have that  $\eE[W_i]=0$ and
$
\sum^n_{i=1}Var(W_i)=n \tau(1-\tau)(f(0))^{-2}\eu^t \eS^{-1}_{{\cal A}^0} \eu$.
Then, by CLT for  independent random variable sequences $(W_i)$, we have
\begin{equation}
\label{cv}
\sqrt{n} f(0) \frac{\eu^t (\widehat \eb_{{\cal A}^0}- \eb^0_{{\cal A}^0})}{\sqrt{\tau(1- \tau) (\eu^t \eS^{-1}_{{\cal A}^0}  \eu)}}  \overset{\cal L} {\underset{n \rightarrow \infty}{\longrightarrow}} {\cal N}(0,1).
\end{equation}
Claim (ii) results taking into account the fact that  $\widehat \eb_{{\cal A}^0}- \eb^0_{{\cal A}^0}=\alpha_n \ed_{{\cal A}^0}$ and relations  (\ref{ad}), (\ref{cv}).
\hspace*{\fill}$\blacksquare$ 

\begin{remark}
The cardinal of the set $\mathcal{A}^0$ may depend on  $n$ and  converge to  $\infty$ as $n \rightarrow \infty$.
\end{remark}

\section{Tuning parameter selection}
In this section we propose a criterion of type BIC to select the tuning parameters $\lambda$ and $\gamma$. This criterion will also estimate the set  $\mathcal{A}^0$, defined by (\ref{A0}). We start with introducing some notations.\\
\noindent $\bullet$ ${\cal A}_n$ a some index set  $\subseteq\{ 1, \cdots, d_n\}$, which does not depend on tuning parameters.\\
$\bullet$ $(\lambda,\gamma) \in (0,\infty)^2$ some  tuning parameters, which does not depend on $n$.\\ 
$\bullet$ $\widehat \eb_{{\cal A}_n}(\lambda,\gamma)$ the seamless $L_0$ quantile estimator of $\eb_{{\cal A}_n}$ obtained on some  index set ${\cal A}_n \subset \{1, \cdots, d_n  \}$ and with $\lambda,\gamma$ as  tuning parameters. We denote its components by $\widehat{\beta}_{\mathcal{A}_n,j}(\lambda,\gamma)$, for $j \in \mathcal{A}_n$.\\
$\bullet$ $\widehat \eb(\lambda_n,\gamma_n)$ the seamless $L_0$ quantile estimator of $\eb$ obtained on the index set  $\{ 1, \cdots, d_n\}$, with $(\lambda_n,\gamma_n)$ as tuning parameters. Then $\widehat \eb(\lambda_n,\gamma_n) = \widehat \eb_n$, with $\widehat \eb_n$ obtained by (\ref{bb}). We denote its components by $\widehat{\beta}_j(\lambda_n,\gamma_n)$, for $j \in \{1, \cdots , d_n \}$.\\
$\bullet$  $\widehat {\cal A}_{\widehat \eb(\lambda_n,\gamma_n)}\equiv \{   j \in \{1, \cdots, d_n \}; \quad \widehat \beta_j(\lambda_n,\gamma_n) \neq 0\}$. \\
$\bullet$ $(\lambda_n,\gamma_n)$ is a tuning parameter sequence such that:
$
\lim_{n \rightarrow \infty} \PP [\widehat {\cal A}_{\widehat \eb(\lambda_n,\gamma_n)} ={\cal A}^0]=1.
$\\

In order to define the BIC criterion, let us consider  $(S_n)_{n \geqslant 1}$,   a sequence of real numbers, defined as:
\begin{itemize}
\item if $d_n/\log n=o(1)$, we consider $S_n=1$ for any $n \in \N$;
\item if $d_n/\log n \not =o(1)$, we consider $(S_n)$ a sequence converging to $\infty$ such that $\displaystyle{ \frac{d_n}{S_n \log n} \rightarrow 0
, \frac{\log n}{ n}   |{\cal A}^0|  S_n \rightarrow 0}$.
\end{itemize}

In order to select ${\cal A}_n$,  $\lambda$ et $\gamma$, we consider the following BIC criterion:
 \begin{equation}
 \label{BIC}
 BIC({{\cal A}_n};{(\lambda,\gamma)}) \equiv \log \big( \frac{1}{n}\sum^n_{i=1} \rho(Y_i - \XX_{i,_{{\cal A}_n}}^t  \widehat \eb_{{\cal A}_n}(\lambda,\gamma))\big)+\frac{\log n}{n} S_n  \| \widehat \eb_{{\cal A}_n}(\lambda,\gamma)\|_0,
 \end{equation}
 with $\| \widehat \eb_{{\cal A}_n}(\lambda,\gamma)\|_0=\sum^{d_n}_{j=1}\e1_{\widehat \beta_{{\cal A}_n,j} (\lambda,\gamma) \neq 0}$.
For the tuning parameters $\lambda_n$, $\gamma_n$ and the estimator $\widehat \eb (\lambda_n,\gamma_n)$, let us consider the  value of the  BIC criterion corresponding to (\ref{BIC}):
\[
BIC(\lambda_n,\gamma_n) \equiv \log \big( \frac{1}{n}\sum^n_{i=1} \rho_\tau \big(Y_i - \XX^t_i  \widehat \eb (\lambda_n,\gamma_n) \big) \big) +\frac{\log n}{n} S_n \| \widehat \eb (\lambda_n,\gamma_n) \|_0.
\]
If the conditions of Theorem \ref{theorem 1} are satisfied, then $ \widehat \eb (\lambda_n,\gamma_n)$ satisfies the  sparsity property,
\[
\lim_{n \rightarrow \infty}\PP[ \{j \in \{1, \cdots, d_n\};\; \widehat \beta_{j} (\lambda_n,\gamma_n) \neq 0 \} = {{\cal A}^0}] = 1 .
\]
In order to prove, by the following theorem, that the  BIC criterion selects correctly, with a probability converging to  1, the tuning parameters $\lambda$ and $\gamma$, we will consider the index sets ${\cal A}_n$ such that $|{\cal A}_n| \leq s_n$, with the assumption $s_n =O(n^a )$, $0 < a < 1/2$. Consider also two index sets   ${A}_{1n}$ et ${A}_{2n}$:
\[{A}_{1n} \equiv  \{{\cal A}_n;  \; {\cal A}^0 \subset {\cal A}_n , {\cal A}^0 \neq {\cal A}_n, |{\cal A}_n| \leq s_n \}, \qquad {A}_{2n} \equiv  \{{\cal A}_n; \; {\cal A}^0  \not \subseteq {\cal A}_n, |{\cal A}_n| \leq s_n \}.
\]
\begin{theorem}
\label{criterion}
We suppose that $0 < \eE[\rho_\tau(\varepsilon)] < \infty$. Then, if instead of assumption (A4) we take $n^{(a-1)/2}d_n^{1/2}\rightarrow 0 $ as $n \rightarrow \infty$, under (A1)-(A3), (A5), we have:
\[
\lim_{n \rightarrow \infty} \PP \big[\min_{\substack{ {\cal A}_n \subseteq \{1, \cdots, d_n\} , (\lambda,\gamma) \in (0,\infty)^2 \\ |{\cal A}_n| \leq s_n }} BIC({{\cal A}_n};{(\lambda,\gamma)})= BIC(\lambda_n,\gamma_n) \big]   =1.
\]
\end{theorem}
\noindent {\bf Proof}.
The theorem is proved if the following two statements are shown:
\begin{equation}
\label{Lemma A1n}
\lim_{n \rightarrow \infty} \PP \big[\min_{{\cal A}_n \in {A}_{1n}} BIC({{\cal A}_n};{(\lambda,\gamma)})> BIC(\lambda_n,\gamma_n) \big]  = 1,
\end{equation}
\begin{equation}
\label{Lemma A2n}
\lim_{n \rightarrow \infty} \PP \big[\min_{{\cal A}_n \in {A}_{2n}} BIC({{\cal A}_n};{(\lambda,\gamma)})> BIC(\lambda_n,\gamma_n) \big]  = 1.
\end{equation}
\underline{\textit{Proof of  relation (\ref{Lemma A1n})}}. 
Since ${\cal A}_n \in {A}_{1n}$, then $|{\cal A}_n | > | {\cal A}^0| $. Let us consider the difference
\[
\begin{array}{l}
 BIC({{\cal A}_n};{(\lambda,\gamma)}) -  BIC(\lambda_n,\gamma_n)   \\
\displaystyle{ =\log \big(1+\frac{n^{-1}\rho_\tau(Y_i - \XX^t_{i, {\cal A}_n}  \widehat \eb_{{\cal A}_n}(\lambda,\gamma)) - n^{-1} \sum^n_{i=1} \rho_\tau(Y_i - \XX^t_{i }  \widehat \eb(\lambda_n,\gamma_n)) }{n^{-1} \sum^n_{i=1} \rho_\tau(Y_i - \XX^t_{i }   \widehat \eb(\lambda_n,\gamma_n))}  \big) + \frac{\log n}{n} S_n \cro{ |{\cal A}_n| - |{\cal A}^0| }}. 
\end{array}
\]
In addition of index set ${\cal A}_n \in {A}_{1n}$, let us consider the following sets: $\widehat{\mathcal{A}}_1 = \{ j; \widehat{\beta}_{n,j} \neq 0 \}$ and  $\widehat{\mathcal{A}}_2 = \{ j; \widehat{\beta}_{\mathcal{A}_n,j} (\lambda,\gamma)\neq 0 \}$. Recall that $\widehat{\eb}_n$ is  $\widehat{\eb}(\lambda_n,\gamma_n)$. Since $\mathcal{A}^0 \subset \mathcal{A}_n$, by Theorem \ref{theorem 1}(i), we have that, $\lim_{n \rightarrow \infty} \PP \cro{\widehat{\mathcal{A}_1}=\widehat{\mathcal{A}_2}=\mathcal{A}^0}=1$. Without loss of generality, we suppose that $\mathcal{A}^0 \subseteq \widehat{\mathcal{A}_1}\subseteq \widehat{\mathcal{A}_2}$, the other cases are similar. Using the elementary inequality $|\rho_\tau(u-v)-\rho_\tau(u) | < |v|$, for all $u,v \in \R$, we have, with probability one,
\begin{eqnarray}
n^{-1} \big| \sum^n_{i=1} \big[ \rho_\tau(Y_i - \XX^t_{i, {\cal A}_n} \widehat \eb_{{\cal A}_n}(\lambda,\gamma)) - \rho_\tau(Y_i - \XX^t_{i }  \widehat \eb(\lambda_n,\gamma_n))\big] \big|  
\leq n^{-1} \sum^n_{i=1} \big| \XX^t_{i, \widehat{\cal A}_2} \pth{\widehat \eb_{{\cal A}_n}(\lambda,\gamma) - \widehat \eb(\lambda_n,\gamma_n)}_{\widehat{\mathcal{A}_2}} \big| \nonumber \\
\leq \max_{1 \leq i \leq n} \|\XX^t_{i, \widehat{\cal A}_2}\|_2   \| \pth{\widehat \eb_{{\cal A}_n}(\lambda,\gamma) - \widehat \eb(\lambda_n,\gamma_n)}_{\widehat{\mathcal{A}_2}} \|_2 \nonumber
\end{eqnarray}
which is, by assumption (A3) and Theorem \ref{Lemma A.1},  $o(\alpha_n) O_{\PP}(\alpha_n)=o_{\PP}(1)$. For the second inequality , the estimators $\widehat \eb_{{\cal A}_n}(\lambda,\gamma)$ were completed by with  zeros for obtaining a vector of dimension $d_n$. Then
\begin{equation}
\label{e2}
n^{-1}\sum^n_{i=1}\rho_\tau(Y_i - \XX^t_{i, {\cal A}_n}  \widehat \eb_{{\cal A}_n}(\lambda,\gamma)) - n^{-1} \sum^n_{i=1} \rho_\tau(Y_i - \XX^t_{i }  \widehat \eb(\lambda_n,\gamma_n)) \overset{\PP} {\underset{n \rightarrow \infty}{\longrightarrow}} 0.
\end{equation}
In the same way, we have:
$
n^{-1} \sum^n_{i=1} \big( \rho_\tau(Y_i - \XX^t_{i }\widehat \eb(\lambda_n,\gamma_n))-\rho_\tau(\varepsilon_i) \big)
\overset{\PP} {\underset{n \rightarrow \infty}{\longrightarrow}} 0$.
On the other hand, be the LLN, 
\begin{equation}
\label{e1}
\frac{1}{n} \sum^n_{i=1} \rho_\tau(\varepsilon_i) \overset{\PP} {\underset{n \rightarrow \infty}{\longrightarrow}} \eE[\rho_\tau(\varepsilon)]  \in (0, \infty).
\end{equation}
Taking into account (\ref{e2}) and (\ref{e1}), we can apply the inequality $\log(1+x) \geq -2 |x|$ for all   $|x|<1/2$,
\begin{eqnarray}
\log \big(1+\frac{n^{-1}\rho_\tau(Y_i - \XX^t_{i, {\cal A}_n}  \widehat \eb_{{\cal A}_n}(\lambda,\gamma)) - n^{-1} \sum^n_{i=1} \rho_\tau(Y_i - \XX^t_{i }  \widehat \eb(\lambda_n,\gamma_n)) }{n^{-1} \sum^n_{i=1} \rho_\tau(Y_i - \XX^t_{i }  \widehat \eb(\lambda_n,\gamma_n))}  \big) \nonumber \\
\qquad \qquad \qquad \geq -2 \frac{\big|n^{-1}\rho_\tau(Y_i - \XX^t_{i, {\cal A}_n}  \widehat \eb_{{\cal A}_n}(\lambda,\gamma)) - n^{-1} \sum^n_{i=1} \rho_\tau(Y_i - \XX^t_{i } \widehat \eb(\lambda_n,\gamma_n)) \big| }{n^{-1} \sum^n_{i=1} \rho_\tau(Y_i - \XX^t_{i }  \widehat \eb(\lambda_n,\gamma_n))}
\label{e3} 
\end{eqnarray}
But, by the proof of Theorem \ref{Lemma A.1}, relation (\ref{Gn}), we have, with probability tending to 1:
$
n^{-1}\rho_\tau(Y_i - \XX^t_{i, {\cal A}_n}  \widehat \eb_{{\cal A}_n}(\lambda,\gamma)) - n^{-1} \sum^n_{i=1} \rho_\tau(Y_i - \XX^t_{i }  \widehat \eb(\lambda_n,\gamma_n))= C\alpha^2_n$,
with $C >0$, for $n$ large enough. Using (\ref{e3}), we have:
\[
\min_{{\cal A}_n \in {A}_{1n}} \big( BIC({{\cal A}_n};{(\lambda,\gamma)})- BIC(\lambda_n,\gamma_n) \big)\geq 
\min_{{\cal A}_n \in {A}_{1n}} \big( -\tilde C \frac{d_n}{n} |{\cal A}_n| +\frac{\log n}{n} C_n \pth{|{\cal A}_n| - | {\cal A}^0|  }\big) > C>0,
\]
with $\tilde C>0$. So, relation (\ref{Lemma A1n}) is proved. \\

\underline{\textit{Proof relation (\ref{Lemma A2n})}}. 
Let be the index sets  ${\cal A}_n \in {A}_{2n}$ and $\widetilde {\cal A}_n \equiv {\cal A}_n \cup {\cal A}^0$.\\
Let $\widehat \eb_{{\cal A}_n}(\lambda,\gamma)$ be the estimator of dimension $|{\cal A}_n|$  built on the   variables $\XX_{{\cal A}_n}$. Let also  $\widetilde \eb_{\widetilde {\cal A}_n}$ equal to $\widehat \eb_{{\cal A}_n}(\lambda,\gamma)$ on ${\cal A}_n$ and completed with 0 to obtain a  vector of dimension $|\widetilde {\cal A}_n|$. Then, denoting $b^0 \equiv \min_{j \in {\cal A}^0} |\beta^0_j| >0 $, we have
\[
\| \widetilde \eb_{\widetilde {\cal A}_n} - \eb^0_{\widetilde {\cal A}_n} \|_2 > \| \eb^0_{{\cal A}^0 \setminus {\cal A}_n} \|_2 \geq b^0 .
\]
Then, since $\rho_\tau(.)$ is convex, we have that there exists $\bar \eb_{\widetilde {\cal A}_n} \in \R^{|\widetilde {\cal A}_n|}$, with $\|\bar \eb_{\widetilde {\cal A}_n} -\eb^0_{\widetilde {\cal A}_n}  \|_2=b^0$, such that
$
\sum^n_{i=1} \rho_\tau(Y_i - \XX^t_{i,{\cal A}_n} \widehat \eb_{{\cal A}_n}(\lambda,\gamma)) \geq \sum^n_{i=1} \rho_\tau(Y_i - \XX^t_{i,\widetilde {\cal A}_n }  \bar \eb_{\widetilde {\cal A}_n})$. 
Thus
\begin{eqnarray}
{\cal R} & \equiv & n^{-1}\sum^n_{i=1} \rho_\tau(Y_i - \XX^t_{i,{\cal A}_n} \widehat \eb_{{\cal A}_n}(\lambda,\gamma)) - \frac{1}{n} \sum^n_{i=1} \rho_\tau(Y_i - \XX^t_{i,\widetilde {\cal A}_n}   \widehat \eb_{\widetilde {\cal A}_n}) \nonumber \\
 & \geq & n^{-1} \big[ \inf_{(\eb - \ebo)_{\widetilde {\cal A}_n} \in  {{\cal B}_\delta(\widetilde {\cal A}_n)}} \eE \big[ G_{n, \widetilde {\cal A}_n} \big( (\eb-\ebo)_{\widetilde {\cal A}_n} \big) \big] -\sup_{(\eb - \ebo)_{\widetilde {\cal A}_n} \in  {{\cal B}_\delta(\widetilde {\cal A}_n)}} \big| G_{n, \widetilde {\cal A}_n} \big( (\eb-\ebo)_{\widetilde {\cal A}_n} \big) \nonumber \\
& &  - \eE \big[ G_{n, \widetilde {\cal A}_n} \big( (\eb-\ebo)_{\widetilde {\cal A}_n} \big) \big]  \big| - G_{n, \widetilde {\cal A}_n}\big( \widehat \eb_{\widetilde {\cal A}_n} - \eb^0_{\widetilde {\cal A}_n} \big) \big], \nonumber
\end{eqnarray}
with
$
G_{n, \widetilde {\cal A}_n}\big( \widehat \eb_{\widetilde {\cal A}_n} - \eb^0_{\widetilde {\cal A}_n} \big) \equiv \sum^n_{i=1} \big[\rho_\tau(Y_i - \XX^t_{i,\widetilde {\cal A}_n}   \widehat \eb_{\widetilde {\cal A}_n}) - \rho_\tau(\varepsilon_i) \big]$  
and $ G_{n, \widetilde {\cal A}_n} \big( (\eb-\ebo)_{\widetilde {\cal A}_n} \big) $ defined similarly.\\
As for the calculation of relation (\ref{eqq2}),  we have, with probability converging to  1:
\[
 \inf_{(\eb - \ebo)_{\widetilde {\cal A}_n} \in  {{\cal B}_\delta(\widetilde {\cal A}_n)}} \eE \big[ G_{n, \widetilde {\cal A}_n} \big( (\eb-\ebo)_{\widetilde {\cal A}_n} \big) \big] \geq C n (b^0)^2, 
\]
with $C>0$ and
\begin{equation}
\label{e6}
 G_{n, \widetilde {\cal A}_n}\big( \widehat \eb_{\widetilde {\cal A}_n} - \eb^0_{\widetilde {\cal A}_n} \big) =O_{\PP}(d_n).
\end{equation}
By Lemma \ref{Lemma A.3} we have:
\[
\sup_{(\eb - \ebo)_{\widetilde {\cal A}_n} \in  {{\cal B}_\delta(\widetilde {\cal A}_n)}} \big| G_{n, \widetilde {\cal A}_n} \big( (\eb-\ebo)_{\widetilde {\cal A}_n} \big) - \eE \big[ G_{n, \widetilde {\cal A}_n} \big( (\eb-\ebo)_{\widetilde {\cal A}_n} \big) \big]  \big| = O_{\PP}\big( d_n^{1/2} n^{(1+a)/2} \big).
\]
Then,  with probability converging to  1, as $n \rightarrow\infty$, we have
$
{\cal R} > C (b^0)^2 - {d_n^{1/2} n^{(1+a)/2}}n^{-1} -  {d_n}n^{-1}$.
Taking into account the assumption $n^{(a-1)/2} d_n^{1/2} \rightarrow 0$, we have that, for $n$ large enough,  with probability converging to  1,
\begin{equation}
\label{e5}
{\cal R} >C (b^0)^2>c_1 >0.
\end{equation}
Hence
\[
\begin{array}{l}
\displaystyle{ 
\min_{{\cal A}_n \in {A}_{2n}} \big[ BIC({{\cal A}_n};{(\lambda,\gamma)}) - BIC({\widetilde {\cal A}_n};{(\lambda,\gamma)})\big]  = \big( |{\cal A}_n | - |\widetilde{{\cal A}_n }|\big)  \frac{\log n}{n} S_n } \\
\displaystyle{ 
+\min_{{\cal A}_n \in {A}_{2n}}  \log \big(1+  \frac{n^{-1}\sum^n_{i=1} \rho_\tau(Y_i - \XX^t_{i,{\cal A}_n}  \widehat \eb_{{\cal A}_n}(\lambda,\gamma)) - n^{-1} \sum^n_{i=1} \rho_\tau(Y_i - \XX^t_{i,\widetilde {\cal A}_n}   \widehat \eb_{\widetilde {\cal A}_n})}{ n^{-1} \sum^n_{i=1} \rho_\tau(Y_i - \XX^t_{i,\widetilde {\cal A}_n}   \widehat \eb_{\widetilde {\cal A}_n})} \big) }
\end{array}
\]
which is,  with probability converging to  1, using (\ref{e5}):
\begin{equation}
\label{eA.21}
\geq \min_{{\cal A}_n \in {A}_{2n}}  \min \big(\log 2, \frac{c_1}{n^{-1} \sum^n_{i=1} \rho_\tau(Y_i - \XX^t_{i,\widetilde {\cal A}_n} \widehat \eb_{\widetilde {\cal A}_n})} \big) - |{\cal A}^0| \frac{\log n}{n} S_n  >0.
\end{equation}
The last inequality  ($>0$) results from  (\ref{e6}) together with  $\eE[\rho_\tau(\varepsilon)] \in (0,\infty)$.\\
As for relation (\ref{Lemma A1n}), we can prove, with probability tending to 1, for $n \rightarrow \infty $:
\begin{equation}
\label{e7}
BIC({\widetilde{\cal A}_n};{(\lambda,\gamma)}) \geq \min_{\substack{{\cal A}'_n\\
{\cal A}'_n \in A_{1n}(2 s_n)}} BIC({{\cal A}'_n};{(\lambda,\gamma)}) \geq BIC(\lambda_n,\gamma_n) ,
 \end{equation}
 with $ A_{1n}(2 s_n)\equiv \{{\cal A}_n;  {\cal A}^0 \subset {\cal A}_n , {\cal A}^0 \neq {\cal A}_n, |{\cal A}_n| \leq 2s_n \}$, $s_n=O(n^a)$, $a \in (0,1/2)$. 
 Then, with probability tending to 1, using (\ref{e7}) and (\ref{eA.21}), we have
 \[
\begin{array}{l}
\displaystyle{ 
 \min_{{\cal A}_n \in {A}_{2n}} BIC({{\cal A}_n};{(\lambda,\gamma)}) -BIC(\lambda_n,\gamma_n)  } \\
 \displaystyle{ = \min_{{\cal A}_n \in {A}_{2n}} \big[BIC({{\cal A}_n};{(\lambda,\gamma)}) -BIC({\widetilde{\cal A}_n};{(\lambda,\gamma)})+BIC({\widetilde{\cal A}_n};{(\lambda,\gamma)})- BIC({{\cal A}^0};(\lambda_n,\gamma_n))\big] }\\
 \displaystyle{\geq \min_{{\cal A}_n \in {A}_{2n}} \big[BIC({{\cal A}_n};{(\lambda,\gamma)}) -BIC({\widetilde {\cal A}_n};{(\lambda,\gamma)}) \big] >0}
\end{array}
\]
and relation (\ref{Lemma A2n}) is proved. 
\hspace*{\fill}$\blacksquare$ \\

Theorem \ref{criterion} implies that we can choose as  tuning parameters $(\lambda_n,\gamma_n)$ thereby:
\[
(\widehat {\cal A}_n, \widehat \lambda_n, \widehat \gamma_n) \equiv \argmin_{\substack{ {\cal A}_n \subseteq \{1, \cdots, d_n\} , (\lambda,\gamma) \in (0,\infty)^2 \\ |{\cal A}_n| \leq s_n }}  BIC({{\cal A}_n};{(\lambda,\gamma)}),
\]
choosing some $s_n$, such that $s_n=O(n^a)$, $a \in (0,1/2)$. 
Obviously $\widehat {\cal A}_n ={\cal A}^0$ with a probability tending to 1. 
Then, in applications, we must first  fix ${\cal A}_n, \lambda,\gamma$ and calculate:
\[
\widehat \eb_{{\cal A}_n}(\lambda,\gamma) =    \argmin_{\eb \in \R^{|{\cal A}_n|}} Q_n(\eb) =  \argmin_{\eb \in \R^{|{\cal A}_n|}}\big( \frac{1}{2 n} \sum^n_{i=1} \rho_\tau(Y_i -\XX^t_{i, {\cal A}_n}   \eb)+ \frac{\lambda}{\log 2} \sum_{j \in {\cal A}_n} \log \big(\frac{|\beta_j|}{|\beta_j|+\gamma}+1 \big)\big)  .
\]
Afterwards,  we vary ${\cal A}_n, \lambda,\gamma$   on grid and take as $\widehat \lambda_n, \widehat \gamma_n$ and $\widehat {\cal A}_n =\widehat {\cal A}_{\widehat \eb(\widehat \lambda_n, \widehat \gamma_n)}$:
\[
(\widehat {\cal A}_n, \widehat \lambda_n, \widehat \gamma_n)  = \displaystyle{  \argmin_{\substack{ {\cal A}_n \subseteq \{1, \cdots, d_n\} , (\lambda,\gamma) \in (0,\infty)^2 \\ |{\cal A}_n| \leq s_n }} \bigg( \log \big( \frac{1}{n} \sum^n_{i=1} \rho_\tau \big(Y_i- \XX^t_{i,{\cal A}_n}  \widehat \eb_{{\cal A}_n}(\lambda,\gamma)  \big) \big)+\frac{\log n}{n} S_n \| \widehat \eb_{{\cal A}_n}(\lambda,\gamma)\|_0 \bigg) } . 
\]

Then, we estimate simultaneously the best  tuning parameters $\widehat \lambda_n$ and $\widehat \gamma_n$ and the parameters $\eb$ that have components different of 0, such that the  corresponding index set $\widehat {\cal A}_n$ is  equal to ${\cal A}^0$, with probability tending to 1. \\

\begin{remark} Theorem \ref{criterion} is the equivalent of Theorem 2 of \cite{LWL12}, where the seamless $L_0$ penalized likelihood approach is considered,  or of  Theorem 2 of \cite{DHL13},  for seamless $L_0$ penalized LS approach.\\
In \cite{LHP14}, a BIC criterion is proposed to select the significant predictor variables of $\XX$ in an high-dimensional quantile model.
\end{remark}

\begin{remark}Algorithm and  numerical part are a very difficult task, since $G_n(\eb_n)$, defined by (\ref{GnB}), is convex in $\eb_n$ and the penalty ${\cal P}en(\eb_n)=\sum^{d_n}_{j=1}p_{SELO}(\beta_j)$ is concave in $\eb_n$, both being continuous, but not differentiable in $\eb_n$. The same type of problem as ours, but with the process $G_n(\eb_n)$ the likelihood (then differentiable in $\eb_n$), was analyzed by \cite{DHL13}. They propose the coordinate descent algorithm to solve the optimization problem. For the method proposed in present paper,  another work should be conducted on numerical method in order to find  the seamless $L_0$ quantile estimator and   the tuning parameters using the criterion given by  Theorem \ref{criterion}.
\end{remark}
\section{Lemmas}
\begin{lemma}
\label{Lemma gh}
Let be the function $g: \R \rightarrow \R$ defined by  $g(x)= \log (h(x)+1)$, with the function $h:\R \rightarrow \R^*_+$, $\displaystyle{h(x)=\frac{|x|}{|x|+\gamma_n}}$. Then, $\forall x_1, x_2 , C\in \R$ such that $|x_1|, |x_2|  \geq C >0$ and $|x_1| -|x_2|  = o(1) $ we have that there exists $\tilde C >0$ such that: 
$
g(x_2)-g(x_1)=\tilde C \gamma_n (|x_2| -|x_1|)  (-1)^{sgn(x_1x_2)}$.
\end{lemma}
\noindent {\bf Proof.} 
By elementary calculus we have
\[
h(x_2)-h(x_1)=\gamma_n \frac{(|x_2|-|x_1|)(-1)^{sgn(x_1x_2)}}{(|x_1|+\gamma_n)(|x_2|+\gamma_n)}.
\]
Then, taking into account the fact that for  $|x| \simeq 0$ we have $\log (x+1) \simeq x$, the lemma follows.
\hspace*{\fill}$\blacksquare$ \\

Let us consider the following notations:
\[
\begin{array}{rll}
\eth & \equiv &(\eb-\ebo)_{{\cal A}_n} , \\
g_{{\cal A}_n}(\varepsilon_i,\eth) &\equiv& \rho_\tau(\varepsilon_i-\XX^t_{i,{\cal A}_n}  \eth) - \rho_\tau(\varepsilon_i) - \eE \left[\rho_\tau(\varepsilon_i-\XX^t_{i,{\cal A}_n}\eth) - \rho_\tau(\varepsilon_i) \right], \qquad \forall i =1, \cdots, n, \\
  {\cal B}_\delta({\cal A}_n) &\equiv &\{ \eth \in \R^{|{\cal A}_n|}; \; \|\eth \|_2 \leq \delta \}, \qquad \forall \delta >0.
\end{array}
\]

\begin{lemma}
\label{Lemma A.3} Under assumptions (A2), (A3), if $s_n=O(n^a)$, with $a \in (0,1/2)$, then, for any
$   \delta > 0$, we have 
\[
\sup_{\substack{{\cal A}_n\\ | {\cal A}_n| \leq s_n}} \sup_{\eth \in {\cal B}_\delta({\cal A}_n)} \big| \sum^n_{i=1}  g_{{\cal A}_n}(\varepsilon_i, \eth)  \big| = O_{\PP} \big( n^{(1+a)/2} d_n^{1/2} \big).
\]
\end{lemma}
\noindent {\bf Proof}. The proof is similar to that of Lemma A.3 of \cite{LHP14}.  
We consider for $k \geq 1$, $\Theta_n(2^{-k}\delta,{\cal A}_n )$ a grid of points in ${\cal B}_\delta({\cal A}_n)$ such that for any $\eth \in {\cal B}_\delta({\cal A}_n)$ there exists $\eth^{(k)} \in \Theta_n(2^{-k}\delta,{\cal A}_n ) $ such that $\|  \eth - \eth^{(k)} \|_2 \leq \delta/ 2^k $. 
If we denote $M \equiv \max_{ 1 \leq i \leq d_n} \| \XX_i \|_2$, then, for a given  constant $C_1 > 0$, let we consider the natural number:
\[
K_n \equiv  \min \big(k \geq 1; \frac{\delta}{2^k} \leq \frac{C_1}{8 M } n^{-1/2} d_n^{1/2} \big).
\]
Using the fact that for any  $u, v \in \R$: $|\rho_\tau(u-v)-\rho_\tau(u) | < |v|$, then, we have with probability 1: 
\begin{equation}
\label{e4}
\sup_{\eth \in {\cal B}_\delta({\cal A}_n)} \big| \sum^n_{i=1} \big[ g_{{\cal A}_n}(\varepsilon_i,\eth) -g_{{\cal A}_n}(\varepsilon_i,\eth^{(K_n)})  \big]  \big| \leq \frac{C_1}{2} n^{1/2} |{\cal A}_n|^{1/2} d_n^{1/2}.
\end{equation}
Denote 
\begin{equation}
\label{P11}
{\cal P}_1  \equiv \PP \big[ \sup_{\substack{{\cal A}_n \\ | {\cal A}_n| \leq s_n}} \sup_{\eth \in {\cal B}_\delta({\cal A}_n)} \big|\sum^n_{i=1}  g_{{\cal A}_n}(\varepsilon_i,\eth) \big|  \geq C n^{(1+a)/2} d_n^{1/2} \big].
\end{equation}
Inequality (\ref{e4}) implies
\[
{\cal P}_1  \leq  \PP \big[  \sup_{\substack{{\cal A}_n \\ | {\cal A}_n| \leq s_n}} \sup_{\eth \in {\cal B}_\delta({\cal A}_n)} \big| \sum^n_{i=1} g_{{\cal A}_n}(\varepsilon_i,\eth^{(K_n)}) \big| \geq \frac{C}{2} n^{(1+a)/2} d_n^{1/2} \big]  .
\]
On the other hand, for the cardinality  $N_k({\cal A}_n)$ of $\Theta_n(2^{-k}\delta,{\cal A}_n )$ , we have:
$
N_k({\cal A}_n) \leq (1 +4 \cdot 2^k)^{|{\cal A}_n|}$. 
Then
\begin{eqnarray}
{\cal P}_1 & \leq & \sum_{{\cal A}_n, |{\cal A}_n| \leq s_n} \PP \big[ \sup_{\eth \in {\cal B}_\delta({\cal A}_n)} \sum^{K_n}_{k=1} |g_{{\cal A}_n}(\varepsilon_i,\eth^{(k)}) - g_{{\cal A}_n}(\varepsilon_i,\eth^{(k-1)}) |  >\frac{C}{2}n^{(1+a)/2}  d_n^{1/2}   \big]\nonumber \\
& \leq &\sum_{{\cal A}_n, |{\cal A}_n| \leq s_n} \sum^{K_n}_{k=1} N_k({\cal A}_n) N_{k-1}({\cal A}_n) {\max}^* \PP \big[ | \sum^n_{i=1} [ |g_{{\cal A}_n}(\varepsilon_i,\eth^{(k)}) - g_{{\cal A}_n}(\varepsilon_i,\eth^{(k-1)}) ] |  \geq \frac{C}{2} \eta_k n^{(1+a)/2}  d_n^{1/2} \big] \nonumber ,
\end{eqnarray}
with $\eta_k > 0$, $\sum^{K_n}_{k=1} \eta_k \leq 1$. The ${\max}^*$ is calculated over all  $\eth^{(k)} \in \Theta_n(2^{-k}\delta,{\cal A}_n )$ and $\eth^{(k-1)} \in \Theta_n(2^{-k+1}\delta,{\cal A}_n )$, with $\|\eth^{(k)} - \eth^{(k-1)} \|_2 \leq 3 \cdot 2^{-k} \delta$. 
Moreover, by assumption (A2) we have:
$
n^{-1} \sum^n_{i=1} \big| \XX_{i,{\cal A}_n}^t  \big(\eth^{(k)} -\eth^{(k-1)} \big)  \big|^2 \leq 18 R_0 2^{-2k} \delta^2$.
We take
$
\eta_k = \max \big( {2^{-k} k^{1/2}}/{8} , {96 R_0^{1/2} 2^{-k} s^{1/2}_n n^{-a/2} d^{-1/2}_n C^{-1} \log^{1/2}(1+4\cdot 2^k)}  \big).
$ 
By the Hoeffding inequality, we obtain:
 \begin{eqnarray}
 {\cal P}_1 & \leq & \sum_{{\cal A}_n; |{\cal A}_n| \leq s_n} \sum^{K_n}_{k=1} \exp \big( 2 s_n \log ( 1+4 \cdot 2^k) - \frac{C^2 \eta^2_k n^a d_n}{48^2 \cdot R_0 \cdot 2^{-2k} \delta^2} \big)   \leq  2\sum_{{\cal A}_n; |{\cal A}_n| \leq s_n} \sum^{K_n}_{k=1} \exp \big( - \frac{C^2 k n^a d_n}{ 2 \cdot 8^2 \cdot 48^2\cdot  R_0 \cdot \delta^2}\big) \nonumber \\
 & \leq &4 \exp \big(s_n \log d_n -  \frac{C^2 k n^a d_n}{ 2 \cdot 8^2 \cdot 48^2 \cdot R_0 \cdot \delta^2} \big)   \rightarrow  0 ,
 \label{P12}
 \end{eqnarray}
 the last  relation following from the fact that $s_n=O(n^a)$. Lemma follows from relations (\ref{P11}) and (\ref{P12}).
\hspace*{\fill}$\blacksquare$ \\

%\textbf{\textsl{References}}

\end{document}